\documentclass[a4paper,11pt]{article}
\pdfoutput=1
\usepackage{a4wide}
\usepackage[english]{babel}
\usepackage{amssymb,amsmath,amsthm}
\usepackage{graphicx}
\usepackage[usenames,dvipsnames]{color}
\usepackage{newcent}
\usepackage[ansinew]{inputenc}
\usepackage[T1]{fontenc}
\usepackage{fouriernc}
\usepackage[pdftex]{hyperref}

\usepackage{type1cm}
\usepackage{eso-pic}

\newcommand{\captionfonts}{\small}

\makeatletter  
\long\def\@makecaption#1#2{%
  \vskip\abovecaptionskip
  \sbox\@tempboxa{{\captionfonts #1: #2}}%
  \ifdim \wd\@tempboxa >\hsize
    {\captionfonts #1: #2\par}
  \else
    \hbox to\hsize{\hfil\box\@tempboxa\hfil}%
  \fi
  \vskip\belowcaptionskip}
\makeatother   

\hypersetup{pdfauthor={C. H. van Dorp}, pdftitle={A free module of Siegel
modular forms}, colorlinks}

\theoremstyle{plain}
\newtheorem{theorem}{Theorem}[section]
\newtheorem{proposition}[theorem]{Proposition}

\newtheorem*{mainthm}{Main Theorem}

\theoremstyle{definition}
\newtheorem{definition}[theorem]{Definition}
\newtheorem{example}[theorem]{Example}

\theoremstyle{remark}
\newtheorem*{remark}{Remark}

\renewcommand{\r}{\mathbf{r}}
\newcommand{\ts}{\!\! \times \!\!}
\newcommand{\PHM}{\mathbf{M}}

\renewcommand{\H}{\mathcal{H}}

\newcommand{\Z}{\mathbb{Z}}
\newcommand{\C}{\mathbb{C}}

\newcommand{\Q}{\mathbb{Q}}
\renewcommand{\epsilon}{\varepsilon}
\renewcommand{\phi}{\varphi}

\newcommand{\Sym}{{\rm Sym}}
\newcommand{\GL}{{\sf GL}}
\newcommand{\SL}{{\sf SL}}
\newcommand{\Sp}{{\sf Sp}}
\newcommand{\pd}{\partial}

\renewcommand{\atop}[2]{\genfrac{}{}{0pt}{1}{#1}{#2}}
\newcommand{\fchoose}[2]{\left(\atop{#1}{#2}\right)}
\newcommand{\mat}[4]{\Bigl(\genfrac{}{}{0pt}{1}{#1}{#3}\,
\genfrac{}{}{0pt}{1}{#2}{#4}\Bigr)}
\newcommand{\Sgroup}{\mathsf{S}}
\newcommand{\gl}{\mathfrak{gl}}
\newcommand{\St}{{\rm St}}
\newcommand{\DRC}{\mathcal{D}}

\title{\sc Generators for a module of  vector-valued  Siegel~modular forms of
degree $2$} 
\author{Christiaan van Dorp\footnote{Korteweg- de Vries Institute for
Mathematics, Universiteit van Amsterdam, Amsterdam. Current address: 
Theoretical Biology and Bioinformatics, Universiteit Utrecht, Utrecht.}}

\begin{document}
\maketitle

\begin{sloppy}

\begin{abstract}\noindent In this paper we will describe all vector-valued
Siegel modular forms of degree $2$ and weight $\Sym^6(\St) \otimes
\det^{k}(\St)$ with $k$ odd. 
These vector-valued forms constitute a module over the ring of classical Siegel
modular forms of degree $2$ and even weight and this module turns out to be
free. 
In order to find generators, we generalize certain Rankin-Cohen differential
operators on triples of classical Siegel modular forms that were first
considered by Ibukiyama and we find a Rankin-Cohen bracket on vector-valued
Siegel modular forms.
\end{abstract}

\section{Introduction}
In comparison to elliptic modular forms, Siegel modular forms and especially
vector-valued Siegel modular forms are much less understood. 
Although the dimensions of the vector spaces of genus $2$ modular forms are
known due to Tsushima \cite{tsushima}, explicit generators are unknown in almost
all cases. 
For some instances however, such generators are known. 
These `first few' examples are due to Satoh \cite{satoh} and Ibukiyama
\cite{ibukiyama1,ibukiyama2}. 

Tsushima's dimension formula can give clues---apart from the dimensions---about
the structure of the modules of forms. 
For instance, the dimension formula can predict relations and in those cases the
modules will not be free. When no relations are predicted, the modules could be
freely generated over the ring of classical modular forms  of even weight.
Ibukiyama has conjectured that one of the modules he studied is in fact free,
but he did not prove this \cite{ibukiyama2}. 
In this paper we will prove his hypothesis and in order to do this, we develop
some methods for constructing vector-valued Siegel modular forms of genus $2$.
Our method allows us to compute a few eigenvalues for the Hecke operators. Our
results agree with calculations done by van der Geer based on his joint work
with Faber \cite{fabergeer1,fabergeer2}.

\subsection{Siegel modular forms}
We will first introduce some notions and notation. Let~$V$ be a finite
dimensional~$\C$-vector space,~$g$ a positive integer and let $\rho : \GL(g,\C)
\rightarrow \GL(V)$ be a representation. 
A Siegel modular form~$f$ of weight~$\rho$ and genus (degree)~$g\geq 2$ is then
a $V$-valued holomorphic function on the Siegel upper half-space~$\H_g$ that
satisfies for every element~$\gamma$ of the symplectic
group~$\Gamma_g:=\Sp(2g,\Z)$ the functional equation
\[
f( \gamma \cdot \tau) = \rho(j(\gamma,\tau)) f(\tau)\, ,\qquad \tau \in \H_g\, .
\]
Here we write $j(\gamma,\tau) := c\tau + d$ for the factor of automorphy and
$\gamma \cdot \tau := (a\tau + b)(c\tau + d)^{-1}$, where $a,b,c,d \in
\gl(g,\Z)$ and $\gamma = \mat{a}{b}{c}{d} \in \Gamma_g$. 
A Siegel modular form $f$ has  a Fourier series 
\[
f(\tau) = \sum_{n} a(n) q^n\, ,\qquad q^n:=e^{2\pi i \sigma(n\tau)}\, ,
\]
where the sum is taken over the set $\Sgroup_g := \left\{\left. n = (n_{ij})\
\right|\ n_{ii} \in \Z, 2n_{ij} \in \Z, n_{ij} = n_{ji} \right\}$ of all
half-integer, symmetric $g\times g$ matrices $n$ and $\sigma(x)$ denotes the
trace of a square matrix $x$. We will write $x'$ for the transpose of a matrix
$x$ and when $y$ is a square matrix of appropriate size, then $y[x] := x'yx$.

The Fourier transform $a_f$ of a Siegel modular form $f = \sum a_f(n)q^n$
defines a function $a_f : \Sgroup_g \rightarrow V$ and if $g>1$, then $a_f(n)
\neq 0 \implies n \succeq 0$ (the Koecher principle). 
Denote by $\Sgroup_g^+$ the subset of $\Sgroup_g$ of semi-positive matrices. A
form $f$ for which $a_f(n)$ vanishes for non-positive $n \in \Sgroup_g^{+}$ is
called a cusp form. 
We will denote the space of modular forms of weight $\rho$ by
$M_{\rho}(\Gamma_g)$ and the space of cusp forms by $S_{\rho}(\Gamma_g)$. 
The general linear group $\GL(g,\Z) \hookrightarrow \Gamma_g$, embedded by $u
\mapsto \mat{u}{0}{0}{{u'}^{-1}}$, acts on~$\Sgroup_g^+$ by $u : n \mapsto unu'$
and the Fourier transform $a_f$ of $f$ behaves well under this action:
\begin{equation}\label{GL2CactionFC}
a_f(unu') = \rho(u)a_f(n)\, .
\end{equation}
From now on, we assume that $g = 2$ unless otherwise specified. 
For convenience, we often write square matrices $\mat{a}{b}{c}{d}$ as
$(a,b;c,d)$ and when no confusion can be possible, we will write $(n,r/2;r/2,m)
\in \Sgroup_2$ as $(n,m,r)$.
The irreducible representations $\rho$ of $\GL(2,\C)$ can be characterized by
their highest weight vector $(\lambda_1 \geq \lambda_2)$ and if $\lambda_2 <0$,
then $\dim_{\C} M_{\rho}(\Gamma_2) = 0$. 
Therefore we only have to consider `polynomial' representations. 
If we write $j = \lambda_1 - \lambda_2$ and $k = \lambda_2$, then $\rho$ will be
isomorphic to $\Sym^j(\St) \otimes \det^k(\St)$, where $\St$ denotes the
standard representation of $\GL(2,\C)$ and $\Sym^j$ and $\det^k$ denote the
$j$-fold symmetric product and $k$-th power of the determinant respectively. 
The representation $\Sym^j(\St)\otimes \det^k(\St)$ is abbreviated to $(j,k)$
and we write~$k$ to denote the weight~$(0,k)$. 

The structure of the ring of classical Siegel modular forms
$M_{\ast}:=\bigoplus_k M_{k}(\Gamma_2)$ of genus~$2$ was determined by Igusa
\cite{igusa1}. 
The subring $M_{\ast}^0 := \bigoplus_{k\equiv 0(2)} M_{k}(\Gamma_2) \subseteq
M_{\ast}$ is a polynomial ring:
\[
M_{\ast}^0 = \C[\phi_4,\phi_6,\chi_{10},\chi_{12}]\, ,
\]
where~$\phi_4$ and~$\phi_6$ are Eisenstein series of weight~$4$ and~$6$ (with
Fourier series normalized at~$(0,0,0)$) and~$\chi_{10}$ and~$\chi_{12}$ are cusp
forms of weight~$10$ and~$12$ (with Fourier series normalized at~$(1,1,1)$). 
Satoh has determined the structure of the $M_{\ast}^0$-module
$M_{(j,\ast)}^i := \bigoplus_{k\equiv i(2)} M_{(j,k)}(\Gamma_2)$ for~$(j,i) =
(2,0)$ and Ibukiyama did the same for~$(j,i) = (2,1)$, $(4,0)$, $(4,1)$ and~$(6,0)$.
In this paper we will determine the structure of~$M_{(6,\ast)}^1$. Our main
result can be formulated as follows.

\begin{mainthm}
The $M_{\ast}^0$-module $M^1_{(6,\ast)}$ is freely generated by seven
elements~$F_{k}$ of weight~$(6,k)$ with $k \in
\left\{11,13,15,17,19,21,23\right\}$.
\end{mainthm}
We shall give the elements mentioned in the above theorem explicitly, but in
order to do this we need to introduce Rankin-Cohen differential operators on
Siegel modular forms. 

\subsection{Rankin-Cohen operators}\label{RCops}
Rankin-Cohen operators (RC-operators) send $t$-tuples of classical Siegel modular forms to
(possibly vector-valued) Siegel modular forms. 
They are therefore useful when we want to find generators for modules of Siegel
modular forms. 
RC-operators were studied in full generality by (among others) Ibukiyama,
Eholzer and Choie \cite{ibukiyama99,eholzeribukiyama,choieeholzer}. 
The general construction of RC-operators can be quite cumbersome and therefore
we only explain the genus $2$ case here. 

Write $H_j = \left\{p \in \C[x,y] \mid \forall\lambda\in \C: p(\lambda x,\lambda
y) = \lambda^j p(x,y)\right\}$ for the space of homogeneous polynomials of
degree $j$ in two variables $x$ and $y$.
The representation $\rho := \Sym^j(\St)\otimes\det^{\ell}(\St) : \GL(2,\C)
\rightarrow \GL(H_j)$ is given explicitly by:
\[
(\rho(G) \cdot p)(x,y) = \det(G)^{\ell} \cdot p((x,y)G)\, , \quad G \in
\GL(2,\C)\, , \quad p \in H_j\, .
\]
Let $R_t = \C[\r_{ij}^s | 0<i\leq j \leq 2,\, 0< s \leq t]$ be the polynomial ring
in the $3t$ variables $\r_{ij}^s$ and consider an element $P$ of
$H_j\otimes_{\C} R_t$. 
It is convenient to write $P(\r_{11}^1,\dots,\r_{22}^t,x,y) =:
P(\r^1,\dots,\r^t;v)$, where $\r^s =
\left(\r_{11}^s,\r_{12}^s;\r_{12}^s,\r_{22}^s\right)$ is a symmetric $2\times 2$
matrix and $v = (x,y)'$. Using this notation, we can give the following
definition.
\begin{definition}
An element $P$ of $H_j \otimes_{\C} R_t$ is called {\em $\rho$-homogeneous} if 
\[
P(G\r^1G',\dots,G\r^t G'; v) = \det(G)^{\ell} \cdot P(\r^1,\dots,\r^t; G'v)
\]
for all $G \in \GL(2,\C)$.
\end{definition}
Ibukiyama refers to certain special elements $P$ of the space
$H_j\otimes_{\C}R_t$ as `associated polynomials', since they are `associated' to
other polynomials $\tilde{P}$ called `pluri-harmonic' polynomials
\cite{ibukiyama99}. 
These polynomials $\tilde{P}$ can be constructed from $P$ by first choosing an
element $k = (k_1,\dots,k_t) \in \Z_{> 0}^t$, which we will refer to as a {\em
type}, and then replacing the matrices $\r^s$ by $\xi^s {\xi^s}'$ with $\xi^s =
(\xi^s_{ij})$ a $2 \times 2k_s$ matrix of indeterminates:
\[
P \mapsto \tilde{P}\ :\  H_j\otimes_{\C} R_t \rightarrow H_j \otimes_{\C}
\C[\xi_{ij}^s \mid 0<i\leq 2, 0< j<2k_s,0<s\leq t]\, ,
\]
\[
\tilde{P}(\xi_{11}^1,\dots,\xi_{2,2k_t}^t,x,y) := P(\xi^1 {\xi^1}',\dots,\xi^t
{\xi^t}';v)\, .
\]
Note that the map $P\mapsto\tilde{P}$ depends on the choice of the type $k$ and that for 
$P \in H_j\otimes_{\C}R_t$, the polynomial $\tilde{P}$ is not necessarily pluri-harmonic.
We can now give the following definition. 
\begin{definition}
An element $P \in H_j \otimes_{\C} R_t$ is called {\em $k$-harmonic} if
$\tilde{P}$ is harmonic in the sense that
\[
\Delta \tilde{P} := \sum_{i,j,s} \frac{\partial^2\tilde{P}}{(\partial
\xi_{ij}^s)^2} =0\, .
\]
\end{definition}
Polynomials that are $\rho$-homogeneous and $k$-harmonic can be used to define
RC-operators as shown in the following theorem and therefore we will refer to
these polynomials as {\em RC-polynomials}. 
The space of RC-polynomials is denoted by $\H_{\rho}(k) \subset H_j \otimes_{\C}
R_t$. 
Write $\tau = (\tau_1,z;z,\tau_2)$ for an element $\tau \in \H_2$ and write $|k|
= \sum_s k_s$.

\begin{theorem}[Ibukiyama]\label{theoremRCops}
Suppose that $P \in \H_{\rho}(k)$ with $\rho = (j,\ell)$ and let $f_1,\dots,
f_t$ be classical Siegel modular forms on $\Gamma_2$ of weight $k_1,\dots,k_t$
respectively. Write $d/d\tau^s := \left(\pd/\pd \tau_{1}^s,\tfrac12\pd/\pd
z^s;\tfrac12 \pd/\pd z^s, \pd/\pd \tau_{2}^s\right)$. 
The $H_j$-valued function
\begin{equation}
\DRC[P](f_1,\dots,f_t)(\tau) := \frac{1}{(2\pi i)^{j/2 +
\ell}}\left.P(d/d\tau^1,\dots,d/d\tau^t; v) f_1(\tau^1)\cdots f_t(\tau^t)
\right|_{\tau^1 = \cdots = \tau^t = \tau}\, ,\quad \tau \in \H_2
\label{defDRC}
\end{equation}
is a Siegel modular form of weight $\rho\otimes \det^{|k|}(\St) = (j,\ell+|k|)$
and genus $2$.
\end{theorem}

A more general version of Theorem~\ref{theoremRCops} was proven by Ibukiyama
\cite{ibukiyama99}. 
Ibukiyama and others also give explicit examples of RC-polynomials and all
RC-polynomials for type of length $2$ were determined explicitly by Miyawaki
\cite{miyawaki}. 
If the length $t$ of the type $k$ equals $2$, then a non-zero Siegel modular
form that has been constructed using a RC-operator will only have weight
$(j,\ell)$ with $\ell$ odd if a classical Siegel modular form of odd weight has
been used in this construction~\cite{satoh,ibukiyama1}. 
The cusp form $\chi_{35} \in M_{35}(\Gamma_2)$ is such a classical Siegel
modular form of odd weight, but if we were to use it in a construction with an
RC-operator,  the weight of the resulting Siegel modular form will have
$\ell\geq 39$. 
The dimension of e.g. $M_{(2,21)}(\Gamma_2)$ equals $1$, hence we need
RC-operators with a higher type length in order to get forms of `low' weight. 
This `trade-off' causes further problems when we increase $j$ and we will give a
(partial) solution to this problem below by introducing a Rankin-Cohen bracket
on vector-valued Siegel modular forms. 

Ibukiyama constructed RC-polynomials of weight $(2,1)$ and $(4,1)$ and type
length $t = 3$ in order to find generators for $M_{(2,\ast)}^1$ and
$M_{(4,\ast)}^1$. We will generalize these polynomials to find generators for
$M_{(6,\ast)}^1$, but our generalization can also be used to find other
vector-valued Siegel modular forms of weight $(j,\ell)$ with $j\geq 2$ and
$\ell\geq 15$ odd.

\subsection{Acknowledgements}
We would like to thank Gerard van der Geer, who stimulated the author to write
this paper, and Fabien Cl\'{e}ry for the useful discussions and their corrections and comments.

\section{Results}
The Hilbert-Poincar\'e series for the dimensions of $M_{(6,\ell)}(\Gamma_2)$
with $\ell$ odd is given by (cf. \cite{tsushima})
\begin{equation}
\sum_{\ell\equiv 1(2)} \dim_{\C} M_{(6,\ell)}(\Gamma_2) \cdot X^{\ell} =
\frac{X^{11} + X^{13} + X^{15} + X^{17} + X^{19} + X^{21} +
X^{23}}{(1-X^4)(1-X^6)(1-X^{10})(1-X^{12})}\,.
\label{dimensionformula}
\end{equation}
This suggests that we should look for generators for $M_{(6,\ast)}^1$ of weights
$(6,\ell)$ with $\ell = 11,13,\dots,23$. 
We will first give two methods for constructing forms of weight $(6,\ell)$ with
$\ell = 15,17,\dots,23$ and $\ell = 11,13$ respectively and then we will give
explicit generators. 

\subsection{RC-polynomials of weight $(j,1)$}

We start with an example by Ibukiyama and Eholzer \cite{eholzeribukiyama}.
The RC-polynomials for elliptic modular forms were studied by Rankin and Cohen
\cite{cohen75} and they can be used to construct RC-polynomials of weight
$(j,0)$. 
Such a genus $1$ RC-polynomial $p_{j,k} \in \C[r_1,r_2]$ can be written as
\begin{equation}\label{cohenspols}
p_{j,k}(r_1,r_2) = \sum_{i=0}^{j/2} (-1)^i\fchoose{j/2}{i} (k_1 + j/2-1)_i (k_2
+j/2-1)_{j/2-i} r_1^{j/2-i} r_2^{i}\, , \qquad j\equiv 0(2)\, .
\end{equation}
Here we use the Pochhammer symbol $(x)_n := x(x-1)(x-2) \cdots (x-n+1)$.
The corresponding RC-polynomial $P_{j,k} \in \H_{(j,0)}(k)$ is then given by
\[
P_{j,k}(\r^1,\r^2;v) := p_{j,k}(\r^1[v],\r^2[v])\,,
\]
and it is easy to verify that~$P_{j,k}$ is indeed $(j,0)$-homogeneous and
$k$-harmonic. 
Given a polynomial~$p$ in the variables $r_1,\dots, r_t$, we then denote the map
that replaces the variables $r_s$ by $\r^s[v]$ by $\Psi$:
\[
\Psi : p(r_1,\dots,r_t) \mapsto p(\r^1[v],\dots,\r^t[v]) : \C[r_1,\dots,r_t]
\longrightarrow \bigoplus_{j\geq 0} H_j \otimes_{\C} R_t\, .
\]
We will now give a construction that uses elliptic (i.e. genus $g=1$)
RC-polynomials and produces RC-polynomials of weight $(j,1)$ and type $k =
(k_1,k_2,k_3)$.
We first consider the cross product on the space of $2 \times 2$ symmetric
matrices. Let $J = (0,1;-1,0)$ and define for $A = (a,b;b,c)$ and $B =
(a',b';b',c')$
\begin{equation}\label{crossprod}
A \ts  B := A J B - B J A = \left(\begin{array}{cc}
2a b' - 2b a' &
a c' - c a' \\
a c' - c a' & 
2b c' - 2c b'
\end{array}\right)\, , 
\end{equation}
then $(G A G')\ts(G B G') = \det(G)\cdot G(A \ts B)G'$ for all $G \in
\GL(2,\C)$. 
\begin{definition}\label{definitionM}
The operator $\PHM_k : \C[r_1,r_2,r_3] \longrightarrow
\bigoplus_{j}H_j\otimes_{\C} R_3$ is defined as follows:
\[
\PHM_k p = 
\r^1\ts\r^2[v]\Psi\left( k_3 p + r_3 \tfrac{\pd}{\pd r_3} p\right) -
\r^1\ts \r^3[v] \Psi \left( k_2 p + r_2 \tfrac{\pd}{\pd r_2} p \right) + 
\r^2 \ts \r^3[v] \Psi \left( k_1 p + r_1 \tfrac{\pd}{\pd r_1} p \right)\, .
\]
\end{definition}
\begin{proposition}\label{propositionM}
Let $p = p_{j,(k_1+1,k_2+1)}$ be the elliptic RC-polynomial defined by
equation~(\ref{cohenspols}), then $\PHM_{k} p$ is a RC-polynomial of weight
$(j+2,1)$ and type $k = (k_1,k_2,k_3)$.
\end{proposition}
The proof of Proposition~\ref{propositionM} is elementary but very tedious. 
We will therefore omit the details. The operator $\PHM_k$ can be shown to
`commute' with the Laplacian $\Delta$ at the cost of a shift in the type
$(k_1,k_2,k_3) \mapsto (k_1+1,k_2+1,k_3+1)$. 
The above result then follows immediately. Note that we could also use other
appropriate polynomials $p \in \C[r_1,r_2,r_3]$.

We now have a recipe to construct vector-valued Siegel modular forms of
weight~$(j,\ell)$ with~$\ell$ odd. 
First take $f_s \in M_{k_s}(\Gamma_2)$ for $s=1,2,3$, then choose an elliptic
RC-polynomial $p = p_{j-2,(k_1+1,k_2+1)}$ and define $P = \PHM_{(k_1,k_2,k_3)} p
\in \H_{(j,1)}(k_1,k_2,k_3)$. 
Theorem~\ref{theoremRCops} then tells us that the function
$\DRC[P](f_1,f_2,f_3)$ is an element of $M_{(j,|k|+1)}(\Gamma_2)$. 
This would be of no use if the resulting Siegel modular forms vanish
identically. 
We can show that $\DRC[P](f_1,f_2,f_3)$ is non-vanishing by computing a non-zero
Fourier coefficient. If $f_s = \sum_n a_s(n) q^n$ and $\DRC[P](f_1,f_2,f_3) =
\sum_n b(n) q^n$, then 
\begin{equation}\label{FC_PHM}
b(n) = \sum_{\atop{(n_1,n_2,n_3) \in (\Sgroup_2^+)^3}{n_1+n_2+n_3 = n}}
P(n_1,n_2,n_3;v)\, a_1(n_1)\, a_2(n_2)\, a_3(n_3)\, .
\end{equation}
Although this is a simple formula, the fact that we take a sum over all triples
$(n_1,n_2,n_3) \in (\Sgroup_2^+)^3$ such that $n_1 + n_2 + n_3 = n$ can result
in a computationally difficult problem, since the number of partitions $n_1 +
n_2 + n_3 = n$ grows very fast as the trace $\sigma(n)$ grows.

\begin{example}
The operator $\PHM_k$ defined in Definition~\ref{definitionM} generalises the
polynomials given by Ibukiyama~\cite{ibukiyama1,ibukiyama2}. These polynomials
are given by 
\[
\PHM_k 1 \in \H_{(2,1)}(k)\quad {\rm and} \quad \PHM_k \bigl((k_1+1) r_2 -
(k_2+1) r_1\bigr) \in \H_{(4,1)}(k)
\]
with $k = (k_1,k_2,k_3)$. Ibukiyama uses these polynomials to define
Rankin--Cohen brackets on triples of classical Siegel modular forms and shows
that the resulting vector-valued  Siegel modular forms generate the modules
$M^1_{(2,\ast)}$ and $M^1_{(4,\ast)}$. 
The operator $\DRC[\PHM_k 1]$ can also be described in terms of the cross 
product~(\ref{crossprod}) and Satoh's RC-brackets \cite{satoh}. Satoh uses the space
of symmetric $2\times 2$ complex matrices as a representation space for
$\Sym^2(\St)$ (where the $\GL(2,\C)$ action is given by $G : A \mapsto G A G'$
for $A$ a symmetric matrix and $G \in \GL(2,\C)$) and then defines
\[
[f_1,f_2] := k_1 f_1 \tfrac{d}{d\tau} f_2 - k_2 f_2 \tfrac{d}{d\tau} f_1 \in
M_{(2,k_1+k_2)}(\Gamma_2), \quad f_s \in M_{k_s}(\Gamma_2)\, .
\]
Let $f_s \in M_{k_s}(\Gamma_2)$ for $s=1,2,3$. Then 
\[
F := [f_1,f_2] \times [f_1,f_3] \in M_{(2,2k_1+k_2+k_3+1)}(\Gamma_2)
\]
and  unraveling the definitions easily shows that $F$ is divisible by $f_1$ in
the $M_{\ast}$-module~$M_{(2,\ast)}$. Ibukiyama's RC-brackets $[f_1,f_2,f_3]$
(in \cite{ibukiyama1}) are then given by
$[f_1,f_2,f_3] := c F/f_1$ for some non-zero constant $c$.
\end{example}

\subsection{A Rankin-Cohen bracket on vector-valued Siegel modular forms}
Our aim is to find all Siegel modular forms of weight~$(6,\ell)$ with~$\ell$
odd. Formula~(\ref{dimensionformula}) shows that $\dim_{\C} M_{(6,11)}(\Gamma_2)
= 1$ and the lowest~$\ell$ for which we can use the operator~$\PHM_k$ to
construct a non-zero form of weight~$(6,\ell)$ equals~$15$. Note that by taking
$k=(4,4,4)$ and $p = p_{4,(5,5)}$, we can get a form $\DRC[\PHM_k p](\phi_4,\phi_4,\phi_4)
\in M_{(6,13)}(\Gamma_2)$, but unfortunately this form vanishes identically. This means
that the above described construction with~$\PHM_k$ is not sufficient for our
purpose, unless we would be able to divide by a classical Siegel modular
form, but this appears to be quite difficult without prior knowledge of the
quotient.
Ibukiyama encountered a similar problem when he determined all modular forms of
weight~$(6,\ell)$ with~$\ell$ even and he solved this by using a theta
series~$\Theta_8$ of weight~$(6,8)$ (cf.~\cite{ibukiyama2,geer}) and a
Klingen-Eisenstein series~$E_6$ of weight~$(6,6)$ (cf.~\cite{arakawa}).
These forms can not be constructed directly by means of RC-operators. 
However, as we will point out below, we can use RC-operators to compute their
Fourier coefficients. In this paper we have scaled $\Theta_8$ such that
the Fourier coefficient of $\Theta_8$ at $(1,1,1)$ equals $x^4y^2 + 2x^3y^3 +
x^2y^4$. The form~$E_6$ is scaled such that its Fourier coefficient at $(1,0,0)$
equals $x^6$. 

The forms~$\Theta_8$ and~$E_6$ can be used to construct forms of weight~$(6,13)$
and~$(6,11)$. 
In order to see how this can be done, we must first give a new interpretation of
the operators $\DRC[\PHM_k p]$, where $p$ is an elliptic RC-polynomial. 
Choose $p = p_{j-2,(k_1+1,k_2+1)}$ and let $f_s \in M_{k_s}(\Gamma_2)$
for~$s=1,2,3$. Also let~$q = p_{j,(k_1,k_2)}$ and define~$F = \DRC[\Psi
q](f_1,f_2) \in M_{(j,k_1+k_2)}(\Gamma_2)$. 
We can re-write~$G := \DRC[\PHM_k p](f_1,f_2,f_3)$ in such a way that 
\[
G =c\cdot \bigl\{F, f_3 \bigr\}
\]
for some bilinear form $\bigl\{ \cdot,\cdot\bigr\}$ on $M_{(j,k_1 +
k_2)}(\Gamma_2) \times M_{k_3}(\Gamma_2)$ and a constant~$c \in \C^{\ast}$. 
We will give the exact description of $\bigl\{ \cdot,\cdot\bigr\}$ below, but
let us first state the advantage of this effort. 
We can replace the modular form~$F$ that was defined via a RC-operator by any
other modular form of weight~$(j,k_1+k_2)$. So for instance, we can take~$F =
E_6$ and~$f_3 = \phi_4$ and if we then apply $\bigl\{ \cdot,\cdot\bigr\}$, we
get
\[
\bigl\{ E_6,\phi_4\bigr\} \in M_{(6,11)}(\Gamma_2)\, .
\]
We will now give an explicit formula for $\bigl\{ \cdot,\cdot\bigr\}$.
\begin{definition}\label{definitionRCbrackets}
Suppose that $F \in M_{(j,k)}(\Gamma_2)$ and $\phi \in M_{\ell}(\Gamma_2)$.
Define the determinant $W(\r)$ by
\[
W(\r) := \left|
\begin{array}{ccc}
\r_{11} & \r_{12} & \r_{22} \\ 
y^2 & -xy & x^2 \\
\frac{\pd^2}{\pd x^2} & \frac{\pd^2}{\pd x \pd y} & \frac{\pd^2}{\pd y^2}
\end{array}
\right| : H_j \rightarrow H_j \otimes_{\C} R_1 \, .
\]
We now define the brackets $\bigl\{ \cdot, \cdot \bigr\}$ as follows.
\[
\bigl\{ F,\phi\bigr\} := \frac{1}{(j-1) 2\pi i} \left((k+j/2-1) W\bigl(\tfrac{d\phi}{d\tau}
\bigr) F - \ell \phi W\bigl(\tfrac{d}{d\tau}\bigr) F\right) \, .
\]
\end{definition}
By our considerations above, we then have the following result.
\begin{proposition}
Take $F$ and $\phi$ as in Definition~\ref{definitionRCbrackets}. The form
$\bigl\{F,\phi\bigr\}$ is then a Siegel modular form of weight $(j,k+\ell+1)$.
\end{proposition}

In order to be able to work with the forms $\bigl\{F,\phi\bigr\}$, we need to
know how to compute Fourier coefficients. 
We can easily find a formula similar to (\ref{FC_PHM}). 
Write $\phi = \sum_{n\succeq 0} a(n)q^n$ and $F = \sum_{n \succeq 0} b(n) q^n$ 
and let $\bigl\{F,\phi\bigr\} = \sum_{n\succeq 0} c(n)q^n$, then we get
\begin{equation}
(j-1) c(n) = \sum_{\atop{(n_1,n_2) \in (\Sgroup_2^+)^2}{n_1+n_2 = n}}  (k+j/2-1)
a(n_1)W(n_1)b(n_2) - \ell a(n_1)W(n_2)b(n_2)\, .
\label{FC_RCbrackets}
\end{equation}

\begin{remark}
A modular form that is defined using the brackets $\bigl\{\cdot,\cdot\bigr\}$
will always be a cusp form. Suppose that we use $F \in M_{(j,k)}(\Gamma_2)$ and
$\phi \in M_{\ell}(\Gamma_2)$ in order to construct $\bigl\{F,\phi\bigr\} \in
M_{(j,k+\ell+1)}(\Gamma_2)$, then either $k + \ell +1$ is odd and hence
$\bigl\{F,\phi\bigr\}$ must be a cusp form, or one of the integers $\ell$ or $k$
is odd, implying that $F$ or $\phi$ is a cusp form. 

We can also see directly from Formula~(\ref{FC_RCbrackets}) that the brackets
map to $S(\Gamma_2)$ by considering the Fourier coefficients at singular indices
$n$. Suppose that $n = n_1 + n_2$ where $n, n_1, n_2 \in S_2^{+}$ and~$n$ is
singular, then we can find a $u \in \SL(2,\Z)$ such that $unu' = (\nu,0,0)$,
$un_1u' = (\nu_1,0,0)$ and $un_2u' = (\nu_2,0,0)$. Hence, by
Formula~(\ref{GL2CactionFC}) we can assume without loss of generality that $n$
is of the form $n = (\nu,0,0)$.  The Fourier coefficients $b(n)$ of $F$ for $n =
(\nu,0,0)$ always have the form $\alpha(\nu) \cdot x^j$ for some constant
$\alpha(\nu)$. Therefore, we get for $s = 1,2$ that
\[
W(n_s) b(n_2) = \nu_s \left|
\begin{array}{cc}
-xy & x^2 \\
\frac{\pd^2}{\pd x \pd y} & \frac{\pd^2}{\pd y^2} 
\end{array}
\right| \alpha(\nu_2) x^j = 0 \, .
\]
This shows that $c(n)$ in Formula~(\ref{FC_RCbrackets}) vanishes for singular
$n$.
\end{remark}

\subsection{Generators for $M_{(6,\ast)}^1$}
As promised, we will now give explicit generators for $M_{(6,\ast)}^1$. 
We first use the brackets $\bigl\{ \cdot,\cdot\bigr\}$ to define forms of weight
$(6,11)$ and $(6,13)$:
\[
F_{11}: = \bigl\{ E_6,\phi_4\bigr\}/1152\, , \qquad F_{13} := \bigl\{
\Theta_8,\phi_4\bigr\}/4
\]
and for the remainder of the generators, we use the construction with $\PHM_k$. 
Choose the following elliptic RC-polynomials $p_i$ and types $k_i$:
\[
\begin{array}{lclcl}
p_i & &{\rm polynomial} & & {\rm type\ } k_i\\ \hline
p_{15}(r_1,r_2) & = & (5r_1^2 - 14 r_1 r_2 + 7 r_2^2) / 160 & \quad  & (5,4,5)\\
p_{17}(r_1,r_2) & = & (4r_1^2 - 8r_1 r_2 + 3 r_2^2)/192 & & (5,6,5)\\
p_{19}(r_1,r_2) & = & (22 r_1^2 - 24 r_1 r_2 + 5 r_2^2)/1920 & & (4,10,4)\\
p_{21}(r_1,r_2) & = & (22 r_1^2 -24 r_1 r_2 + 5 r_2^2)/2880 & & (4,10,6)\\
p_{23}(r_1,r_2) & = & (13 r_1^2 - 14 r_1 r_2 + 3 r_2^2)/16 & & (5,12,5)\\
\end{array}
\]
We then write $P_i = \PHM_{k_i} p_i$ for $i = 15,17,\dots,23$ and define modular
forms $F_i$ as follows.
\[
\begin{array}{lclcl}
F_i & & {\rm modular\ form} &\quad & {\rm weight}\\ \hline
F_{15} &=& \DRC[P_{15}](\chi_5,\phi_4,\chi_5) & & (6,15)\\
F_{17} &=& \DRC[P_{17}](\chi_5,\phi_6,\chi_5) & & (6,17)\\
F_{19} &=& \DRC[P_{19}](\phi_4,\chi_{10},\phi_4) & & (6,19)\\
F_{21} &=& \DRC[P_{21}](\phi_4,\chi_{10},\phi_6) & & (6,21)\\
F_{23} &=& \DRC[P_{23}](\chi_5,\chi_{12},\chi_{5}) & & (6,23)
\end{array}
\]
The form $\chi_5$ denotes the square root of $\chi_{10}$ in the ring of
holomorphic functions on $\H_2$ (see e.g.~\cite{freitag}). 
Our main result can now be formulated as follows:
\begin{mainthm}
The $M_{\ast}^0$-module $M_{(6,\ast)}^1$ is free and can be written as the
following direct~sum:
\[
M_{(6,\ast)}^1 = \bigoplus_{i\in I} F_i \cdot M_{\ast}^0\, ,
\]
where $F_i$ are defined above and the direct sum is taken over the set $I =
\{11,13,15,17,19,21,23\}$.
\end{mainthm}

\subsection{Eigenvalues of the Hecke operators}
Since we now know all modular forms of weight $(6,k)$ with $k\in \Z$, we can
calculate the eigenvalues of the Hecke operators $T(p)$ (cf.
\cite{arakawa,andrianov}). 
We did this for $p=2,3$ and some $k$ (Table~\ref{table_eigenvalues2}
and~\ref{table_eigenvalues3}).

\begin{table}[t]
\begin{center}
\begin{tabular}{l|l|l}
$k$ & $\lambda(2)$ on $N_{(6,k)}(\Gamma_2)$ & $\lambda(2)$ on $S_{(6,k)}(\Gamma_2)$\\
\hline
$6$ 	& $-24 \cdot (1+2^4)$ 		& --- 				\\
$8$ 	& --- 						& $0$ 				\\
$10$ 	& $216\cdot(1+2^8)$		& $1680$ 			\\
$11$     & --- 						& $-11616$			\\
$12$	& $-528\cdot(1+2^{10})$	& $X^2 - 22368 X + 57231360$			\\
$13$	& --- 						& $-24000$			\\
$15$  	& --- 						& $X^2 + 68256 X + 593510400$ \\
$17$	& --- 						& $X^3 + 363264X^2 + 136028160 X - 4603543289856000$\\
$19$	&--- 						& $X^4 + 1202400 X^3 -  1311202861056 X^2 $\\ & &\qquad $ - 179858880190218240 X - 1566691549034368204800$\\
\end{tabular}
\end{center}
\caption{Eigenvalues of the Hecke operator~$T(2)$ on $M_{(6,k)}(\Gamma_2)$ for some values of~$k$. If a polynomial in~$X$ is given, the eigenvalues~$\lambda(2)$ are the roots of this polynomial. The space $N_{(6,k)}(\Gamma_2)$ is the orthogonal complement of $S_{(6,k)}(\Gamma_2)$ with respect to the Petersson product.}\label{table_eigenvalues2}
\end{table}

\begin{table}
\begin{center}
\begin{tabular}{l|l}
$k$ & $\lambda(3)$ on $S_{(6,k)}(\Gamma_2)$\\
\hline
$8$ & $-27000$ \\
$10$ & $-6120$ \\
$11$ & $-106488$ \\
$12$ & $X^2 + 335664 X - 14832719455680$\\
$13$ & $-8505000$\\
$15$ & $X^2 + 228022128 X + 8319716602228800$\\
$17$ & $X^3 + 1086146712 X^2 - 341960280255362880 X -188775313801934579676864000$
\end{tabular}
\end{center}
\caption[Table]{Eigenvalues of the Hecke operator $T(3)$ on $M_{(6,k)}(\Gamma_2)$ for some values of $k$.}\label{table_eigenvalues3}
\end{table}

At our request, G. van der Geer computed some of these eigenvalues using a
completely independent method that is based on counting points on hyperelliptic
curves over finite fields~\cite{geer,fabergeer1,fabergeer2}. 
The values listed here agree with these. 
Also note that the characteristic polynomials of $T(2)$ and $T(3)$
on~$S_{(6,k)}(\Gamma_2)$ with $k=12$ and $15$ have the following discriminant:
\[
\begin{array}{l|l|l}
k & \Delta(\det(T(2)-X)) & \Delta(\det(T(3)-X))\\
\hline
12 & 2^{10} 3^2 7^2 601           &    2^{14} 3^6 7^2 13^2 601\\
15 & 2^{10} 3^2 29 \cdot 83 \cdot 103  &   2^{12} 3^8 29\cdot 53^2 83 \cdot 103\\
\end{array}
\]
This shows  for $p=2,3$ that the eigenvalues of $T(p)$ on $S_{(6,12)}(\Gamma_2)$
and $S_{(6,15)}(\Gamma_2)$ are elements of the same quadratic number field
$\Q(\sqrt{601})$ and $\Q(\sqrt{29\cdot 83 \cdot 103})$ respectively. 
We also verified that the characteristic polynomials of $T(2)$ and $T(3)$ on
$S_{(6,17)}(\Gamma_2)$ define the same number field.

\section{Proof of the main theorem}
We will only give a sketch of the full proof since it involves many elementary
computations.
Let~$U$ denote the Hodge bundle corresponding to the factor of automorphy $j$. 
The forms~$F_i$ defined above are sections of $\Sym^6(U) \otimes L^{\otimes
k_i}$, where $L$ denotes the line bundle $\det(U)$. 
In order to show that the forms $F_i$ are independent over $M_{\ast}^0$, we have
to show that $\chi_{140}:=F_{11} \wedge F_{13} \wedge \dots \wedge F_{23} \in
L^{\otimes 140}$ is non-vanishing. 
The form $\chi_{140}$ is an element of $M_{140}(\Gamma_2)$ and we can compute
the Fourier coefficients of this form. 
We will therefore prove the main theorem by computing a non-zero Fourier
coefficient of $\chi_{140}$. 
Write $\chi_{140} = \sum_n c(n) q^n$ and $F_i = \sum_n a_i(n) q^n$.  The
following formula holds for the Fourier coefficients $c(n)$:
\begin{equation}\label{FCchi140}
c(n) = \sum_{\sum_i n_i = n} \det(a_{11}(n_1),a_{13}(n_2),\dots,a_{23}(n_7))\,
. 
\end{equation}
As mentioned above, this can be hard to compute if $\sigma(n)$ is large. In
order to find a non-zero Fourier coefficient $c(n)$ of $\chi_{140}$, we need at
least $\sigma(n) \geq 14$ since the forms $F_i$ are cusp forms. Indeed, under
the Siegel operator $\Phi$, the forms $F_i$ map to elliptic modular forms of odd
weight. These forms all vanish.

Using algorithms provided by Resnikoff and  Salda\~{n}a \cite{resnikoff74}, we
calculated Fourier coefficients of the classical Siegel modular forms
$\phi_4,\phi_6,\chi_{10}$ and $\chi_{12}$. 
The Fourier coefficients of $\chi_5$ can be computed using $\chi_5^2 =
\chi_{10}$.\footnote{While computing the Fourier coefficients of $\chi_5$ we
encountered an error in Table~IV of \cite{resnikoff74}. The coefficients at
calculation classes $(3,3,3)$ and $(2,6,0)$ should have opposite sign.} 
We then were able to compute Fourier coefficients of the forms~$F_i$. 

In order to find Fourier coefficients of $F_{11}$ and $F_{13}$, we first had to
compute Fourier coefficients of~$E_6$ and~$\Theta_8$. 
The modular form $\phi_4 E_6$ is an element of $M_{(6,10)}(\Gamma_2)$ and this
space has dimension~$2$. 
Using a RC-polynomial of weight $(6,0)$ and the modular forms $\phi_4$ and
$\phi_6$, we can find a modular forms $F_{10}$ in the space
$M_{(6,10)}(\Gamma_2)$ (cf. Ibukiyama \cite{ibukiyama2}) and this form is not an
eigenform. 
Hence, we can find Fourier coefficients of $F_{10}$ and $T(2)F_{10}$. For some
$\alpha,\beta \in \mathbb{Q}$, we must have $\alpha F_{10} + \beta T(2) F_{10} = \phi_4
E_6$. The form $E_6$ is an eigenform and by a theorem due to
Arakawa~\cite{arakawa} we know the corresponding eigenvalue of $T(2)$ (cf.
Table~\ref{table_eigenvalues2}). 
This gives a relation for the Fourier coefficients of $E_6$ and using this
relation, we were able to find $\alpha$ and $\beta$. This also allowed us to
compute the Fourier coefficients of $E_6$.

A similar method can be used to compute Fourier coefficients of $\Theta_8$.
Again using a RC-operator of weight $(6,2)$ and the modular forms $\phi_4$ and
$\phi_6$, we can find a modular form $F_{12}$ in $S_{(6,12)}(\Gamma_2)$ (cf.
Ibukiyama \cite{ibukiyama2}). This form is again not an eigenform. The form
$\phi_4\Theta_8 \in S_{(6,12)}(\Gamma_2)$ must be a linear combination $\alpha
F_{12}+\beta T(2)F_{12}$ and since Ibukiyama has computed a few Fourier
coefficients of $\Theta_8$ (cf. \cite{ibukiyamageer}), we were able to determine
$\alpha$ and $\beta$. 

We then used formulas~(\ref{FC_PHM}) and~(\ref{FC_RCbrackets}) to compute
Fourier coefficients of the forms $F_i$. A few examples are given in
Table~\ref{FCtable}.

\begin{table}
\[
\begin{array}{l|ccccccc}
i & 11 & 13 & 15 & 17 & 19 & 21 & 23\\ \hline
n & (1,1,0) & (1,1,1) & (2,1,0) & (2,1,0) & (2,1,1) & (2,1,1) & (2,2,1) \\
\hline
a_i(n) & 0 	&  0     	& 0 	&   0   	&   1   	&  -5   
  &  3\\
& -20 	&  -2	& 312 	&   0  	&  14  	&  -10    &    -37\\
& 0	&  -5    	&  0 	&   0     &   36   	&    -6     &  -50\\
& 0  	&  0     	&  180  & -300 	&   24  	&    -24   &     0\\
& 0  	&  5    	& 0		&   0    	&     0  	&    -30
   &   50\\
& 20  	&  2   	& -102  &  354  	&    0    &    -12    &  37\\
& 0   	&  0      	&   0    	&   0  	&    0   	&    0   	
 &       -3
\end{array}
\]
\caption{\label{FCtable} A few Fourier coefficients of the forms $F_i$. The
Fourier coefficients are written as column vectors. The determinant of the above
matrix occurs in the sum~(\ref{FCchi140}) for $c(12,8,4)$. This determinant
equals $2^{14}3^5 5^3 11$.  }
\end{table}

We wrote a script to compute the Fourier coefficient of $\chi_{140}$ at $n =
(12,8,4)$ and found that $c(12,8,4) = -2^{18}3^75^2 \neq 0$. We checked our
computations by also computing the Fourier coefficient at $(12,8,-4)$ and found
that $c(12,8,-4) = c(12,8,4)$ which is in line with
equation~(\ref{GL2CactionFC}). This proves our main result.

\end{sloppy}
\end{document}